\documentclass[11pt,twoside]{amsart}
\usepackage{amssymb}

\usepackage[latin1]{inputenc}
\usepackage[T1]{fontenc}
\usepackage{mathtools}
\usepackage{graphicx}
\usepackage{tikz}
\usepackage{mathabx}
\usetikzlibrary{chains}

\PassOptionsToPackage{pdfusetitle,pagebackref,colorlinks}{hyperref}
\usepackage{bookmark}
\hypersetup{
  linkcolor={red!70!black},
  citecolor={green!70!black},
  urlcolor={blue!80!black}
}

\usepackage[latin1]{inputenc}
\usepackage{amsmath}
\usepackage{amsthm}

\usepackage[all]{xy}
\usepackage{hyperref}
\newtheorem{thm}{Theorem}[section]

\numberwithin{equation}{section}

\theoremstyle{definition}

\newtheorem{remark}[thm]{Remark}
\newtheorem{conj}[thm]{Conjecture}

\newtheorem{question}[thm]{Question}

\newcommand{\Db}{{\rm D}^{\rm b}}

\newcommand{\Aut}{{\rm Aut}}

\newcommand{\Br}{{\rm Br}}

\newcommand{\Hom}{{\rm Hom}}

\newcommand{\id}{{\rm id}}

\newcommand{\cal}{\mathcal}
\newcommand{\ka}{{\cal A}}

\newcommand{\kc}{{\cal C}}

\newcommand{\ke}{{\cal E}}

\newcommand{\ko}{{\cal O}}

\newcommand{\kt}{{\cal T}}

\newcommand{\ZZ}{\mathbb{Z}}
\newcommand{\QQ}{\mathbb{Q}}

\newcommand{\CC}{\mathbb{C}}

\newcommand{\PP}{\mathbb{P}}

\renewcommand{\to}{\xymatrix@1@=15pt{\ar[r]&}}
\newcommand{\lto}{\xymatrix@1@=15pt{&\ar[l]}}
\renewcommand{\rightarrow}{\xymatrix@1@=15pt{\ar[r]&}}
\renewcommand{\mapsto}{\xymatrix@1@=15pt{\ar@{|->}[r]&}}
\newcommand{\mapslto}{\xymatrix@1@=15pt{&\ar@{|->}[l]&}}
\renewcommand{\twoheadrightarrow}{\xymatrix@1@=18pt{\ar@{->>}[r]&}}
\renewcommand{\hookrightarrow}{\xymatrix@1@=15pt{\ar@{^(->}[r]&}}
\newcommand{\hook}{\xymatrix@1@=15pt{\ar@{^(->}[r]&}}
\newcommand{\congpf}{\xymatrix@1@=15pt{\ar[r]^-\sim&}}
\renewcommand{\cong}{\simeq}
\newcommand{\TBC}[1]{}
\makeatletter
\def\blfootnote{\xdef\@thefnmark{}\@footnotetext}
\makeatother

\begin{document}

\title{The K3 category of a cubic fourfold -- an update}

\author[D.\ Huybrechts]{Daniel Huybrechts}

\address{Mathematisches Institut \& Hausdorff Center for Mathematics,
Universit{\"a}t Bonn, Endenicher Allee 60, 53115 Bonn, Germany}
\email{huybrech@math.uni-bonn.de}

\begin{abstract} \noindent
We revisit  \cite{HuyComp}, review subsequent developments and highlight some open questions. The exposition avoids the more technical points and concentrates
on the main ideas and the overall picture.
 \vspace{-2mm}
\end{abstract}

\maketitle
\blfootnote{These are the notes of a talk
delivered at the prize ceremony of the Compositio Foundation in Berlin  on 10 February 2022. The author is supported by the ERC Synergy Grant HyperK (ID 854361).}\blfootnote{I wish to thank N.\ Addington, P.\ Belmans, and
R.\ Thomas for helpful comments and suggestions.}

\marginpar{}

The paper \cite{HuyComp} relies crucially on the work of Hassett \cite{HassettComp}, Kuznetsov \cite{Kuz1,Kuz3}, Addington--Thomas \cite{AT}, and, going back further, on the work of Mukai \cite{MukaiVB} and Orlov \cite{OrlovK3}. \smallskip

The purpose of \cite{HuyComp} was twofold: (i) Upgrade the existing theory linking cubic fourfolds and K3 surfaces to incorporate twisted K3 surfaces; (ii)
Establish the analogue of the theory of derived equivalences between (twisted) K3 surfaces for Kuznetsov's K3 categories associated with cubic fourfolds. 

The paper \cite{HuyComp} succeeded with regard to (i) but was only partially successful concerning (ii). Even today, as will be reviewed below, the Mukai--Orlov program for cubic fourfolds remains incomplete.

For more comprehensive treatments of the general theory we refer to \cite{HassettSurvey,HuySurvey,HuyCubic,MSSurvey}.
\section{Hodge theory and derived categories}
We recall the basic features of the two main players: Hodge structures of cubic fourfolds and their derived categories.

\subsection{Hodge theory}\label{sec:Hodge}
That Hodge structures of cubic fourfolds show
certain similarities to  Hodge structures of K3 surfaces was first observed by
Deligne and Rapoport \cite{DR} and later further exploited by Beauville and Donagi
\cite{BD}, Voisin \cite{Voisin}, and others. The first systematic study of this mysterious link was undertaken by Hassett in his PhD thesis \cite{HassettComp}.

For a smooth cubic fourfold $X\subset \PP^5$ over the complex numbers the middle
cohomology $H^4(X,\ZZ)\cong\ZZ^{\oplus 23}$ is of rank $23$ and its non-trivial Hodge numbers are $h^{3,1}=h^{1,3}=1$ and $h^{2,2}=21$. Thus, up to Tate twist,
it is a weight two Hodge structure of K3 type. The rank can be reduced to $22$ by
passing to the primitive cohomology $H^4(X,\ZZ)_{\rm pr}\coloneqq {h^2}^\perp\subset H^4(X,\ZZ)$. As a lattice, the latter can be embedded primitively and isometrically (up to a global sign)  into the full cohomology lattice of a K3 surface:
$$H^4(X,\ZZ)_{\rm pr}\,\hookrightarrow H^\ast({\rm K3},\ZZ).$$
Note that the complement $H^4(X,\ZZ)_{\rm pr}^\perp\subset H^\ast(X,\ZZ)$, which is
the rank five lattice $H^0\oplus H^2\oplus \ZZ\cdot h^2\oplus H^6\oplus H^8$,
carries a trivial Hodge structure, which, in particular, is constant under deformations of $X$.

\begin{remark}\label{rem:GT}
The global Torelli theorem for cubic fourfolds proved by Voisin \cite{Voisin} states
that two  smooth cubic fourfolds $X,X'\subset\PP^5$ are isomorphic if and only if
there exists a Hodge isometry $H^4(X,\ZZ)_{\rm pr}\cong H^4(X',\ZZ)_{\rm pr}$:
$$X\cong X'~\Leftrightarrow~H^4(X,\ZZ)_{\rm pr}\cong H^4(X',\ZZ)_{\rm pr}.$$
Also, the group $\Aut(X)$  of automorphisms of $X$ and $\Aut(H^4(X,\ZZ)_{\rm pr})$ of Hodge isometries of  $H^4(X,\ZZ)_{\rm pr}$ coincide up to the map $-\id$.
\end{remark}

The Hodge numbers being $h^{3,1}=h^{1,3}=1$ and the lattice $H^4(X,\ZZ)_{\rm pr}$ closely linked to the cohomology ring of a K3 surface, the following 
question arises naturally.\smallskip

\begin{question}\label{qu:Hassett}
Under what conditions is the Hodge structure
$H^4(X,\ZZ)_{\rm pr}$ realised (in a sense to be made precise) by a K3 surface?\end{question}

Hassett \cite{HassettComp} approached this question via the primitive cohomology of
polarised K3 surfaces $(S,H)$. The weight two Hodge structure $H^2(S,\ZZ)_{\rm pr}$
is of rank $21$ and he asked under which conditions it can embedded, primitively, isometrically, and compatibly with the Hodge structures, into the primitive cohomology of $X$:
\begin{equation}\label{eqn:Hassett}H^2(S,\ZZ)_{\rm pr}\,\hookrightarrow H^4(X,\ZZ)_{\rm pr}
\end{equation}
(which automatically is of corank one).
The approach that turned out to be better adapted to the categorical viewpoint asks
instead for a Hodge isometry
\begin{equation}\label{eqn:Hassettmodern}
\widetilde H(X,\ZZ)\cong \widetilde H(S,\ZZ).
\end{equation}
Here, $\widetilde H(S,\ZZ)$ is the full cohomology $H^\ast(S,\ZZ)$ viewed as a weight two Hodge structure (the Mukai Hodge structure).
The definition of the left hand side is more involved. Addington and Thomas \cite{AT}
introduced a corank three lattice of the topological K-theory
$$
\widetilde H(X,\ZZ)\subset K_{\rm top}(X)$$
as the complement of $[\ko],[\ko(1)],[\ko(2)]\in K_{\rm top}(X)$ and
gave it a weight two Hodge structure   by declaring $\widetilde H^{2,0}(X)$ to be the pull-back (via the Chern character) 
of $H^{3,1}(X)$. The pairing is defined in terms of the Mukai vector, i.e.\ the Riemann--Roch pairing.
The passage from the standard integral cohomology $H^\ast(X,\ZZ)$ to integral
topological K-theory is subtle.
For K3 surfaces the two coincide and observing the importance of the difference was crucial for the development of the theory.
The two sides of a Hodge isometry (\ref{eqn:Hassettmodern}) extend the primitive cohomologies
of $X$ and $S$:
$$H^4(X,\ZZ)_{\rm pr}\subset \widetilde H(X,\ZZ)\text{ and } H^2(S,\ZZ)_{\rm pr}\subset \widetilde H(S,\ZZ).$$
We think of (\ref{eqn:Hassettmodern}) as the unpolarised version of
(\ref{eqn:Hassett}), i.e.\ it does not involve the polarisation of the K3 surface $S$. There is of course no ambiguity in the choice of the polarisation of $X$.

\subsection{Derived categories} Instead of studying the Hodge structure of
a cubic fourfold, one can look at its bounded derived category $\Db(X)$ of coherent
sheaves.

\begin{remark}\label{rem:BO}
The analogue of Remark \ref{rem:GT} is the result by Bondal and Orlov \cite{BO}:
Two smooth cubic fourfolds $X, X'$ are isomorphic if and only if there exists
an exact linear equivalence $\Db(X)\cong \Db(X')$:
$$X\cong X'~\Leftrightarrow ~ \Db(X)\cong \Db(X').$$
Similarly, the groups $\Aut(X)$ of automorphisms  of $X$ and $\Aut(\Db(X))$
of exact linear auto-equivalences  of $\Db(X)$ essentially coincide,
up to shifts and tensor product with line bundles.
\end{remark}

A link to K3 surfaces is established via a certain full triangulated subcategory
$$\ka_X\subset \Db(X).$$ By definition, it is defined by the orthogonality condition
$\Hom(\ko(i),E[\ast])=0$ for $i=0,1,2$. The crucial observation, a celebrated result
of Kuznetsov \cite{Kuz1,Kuz3}, asserts that with this definition $\ka_X$ has a property that is typical
for the derived category $\Db(S)$ of a K3  surface or an abelian surface: For all $E,F\in \ka_X$
there exists a functorial isomorphism $\Hom(E,F)\cong \Hom(F,E[2])^\ast$.
Another interpretation of $\ka_X$ as a category of graded matrix factorisations
was obtained by Buchweitz \cite{Buchweitz} and Orlov \cite{OrlovMF}.

The analogue of Question \ref{qu:Hassett} is then the following.

\begin{question}
Under what conditions on $X$ does there exist an exact linear equivalence $\ka_X\cong \Db(S)$ to the derived category of a K3 surface $S$?
\end{question}

\section{Rationality conjectures} The  middle cohomology
$H^3(Y,\ZZ)$ and the intermediate Jacobian $J(Y)$ of a cubic threefold
$Y\subset\PP^4$ are important and interesting invariants beyond their 
celebrated application to the irrationality of smooth cubic threefolds due to
Clemens and Griffiths. The same is true for the middle cohomology 
$H^4(X,\ZZ)$ of a cubic fourfold and for the K3 category $\ka_X$.\footnote{Although there are some similarities between the case of cubics of dimension three and four, the role of their Kuznetsov components
is different. For example, a smooth cubic threefold is determined by
its Kuznetsov component \cite{BMMS}.}\smallskip

For those among the readers who need a concrete application as a motivation to
study the two notions $H^4(X,\ZZ)$ and $\ka_X$, we recall two popular rationality conjectures
(and two less popular ones). 

\begin{conj}[Hassett]\label{conj:Hassett} A smooth cubic fourfold $X\subset \PP^5$ is rational if and only
if there exists a primitive  isometric 
embedding of Hodge structures $H^2(S,\ZZ)_{\rm pr}\,\hookrightarrow H^4(X,\ZZ)_{\rm pr}$
for some polarised K3 surface $(S,H)$:
\begin{equation}\label{eqnHC}
X\text{ rational } \Leftrightarrow~  \exists \, \, H^2(S,\ZZ)_{\rm pr}\,\hookrightarrow H^4(X,\ZZ)_{\rm pr}.
\end{equation}

Equivalently, $X$ is rational if and only
if there exists a Hodge isometry $\widetilde H(X,\ZZ)\cong \widetilde H(S,\ZZ)$
for some (unpolarised) K3 surface $S$
\begin{equation}\label{conj:AT} X\text{ rational } \Leftrightarrow ~ \exists \, \, \widetilde H(X,\ZZ)\cong \widetilde H(S,\ZZ).
\end{equation}
\end{conj}

Hassett \cite{HassettComp} avoids stating the conjecture (\ref{eqnHC}) explicitly, but it is clearly what he (and Harris) had in mind, and the reformulation (\ref{conj:AT}) is due to Addington and Thomas.   See Remark \ref{rem:Ratlatest} for a collection of known results.\smallskip

The next conjecture is stated explicitly for the first time in \cite[Conj.\ 3]{Kuz5}.

\begin{conj}[Kuznetsov]\label{conj:Kuznetsov} A smooth cubic fourfold $X\subset \PP^5$ is rational
if and only if there exists an exact linear equivalence $\ka_X\cong \Db(S)$ for some
K3 surface $S$:
\begin{equation}
X\text{ rational } \Leftrightarrow ~\exists\,\, \ka_X\cong\Db(S).
\end{equation}
\end{conj}

The two conjectures are of a completely different flavour, but the striking result of Addington and Thomas \cite{AT} shows that they are in fact equivalent. We will comment on this in more detail further below.\smallskip

But of course, until either of the two conjectures is settled (and hence both), we are free 
to propose alternative ones. To the more cautious reader, for whom the two conjectures above predict
rationality for too many cubics, the following version might appeal \cite{GS}.
The geometric evidence for this rationality conjecture seems somewhat
stronger than for the previous two.

\begin{conj}[Galkin--Shinder]\label{conj:GS}
A smooth cubic fourfold $X$ is rational if and only if
its Fano variety of lines $F(X)$ is birational to the Hilbert scheme $S^{[2]}$ of some K3
surface $S$:
\begin{equation}X \text{ rational } \Leftrightarrow~\exists\,\,  F(X)\sim S^{[2]}.
\end{equation}
\end{conj}

If, as conjectured by Hassett and Kuznetsov, K3 surfaces really decide about the
rationality of a cubic fourfold, one could venture a conjecture predicting more cubics to be rational by allowing twisted K3 surfaces $(S,\alpha\in \Br(S))$. This then would lead to the following more provocative speculation.

\begin{conj} A smooth cubic fourfold $X$ is rational if and only if
there exists an exact linear equivalence $\ka_X\cong\Db(S,\alpha)$ for some
twisted K3 surface $(S,\alpha)$.
\end{conj}

The only evidence against this last rationality conjecture is that it would predict all
cubic fourfolds containing a plane to be rational. However, as those are among the most studied examples of cubic fourfolds, someone would have proved their rationality by now. Note that among cubic fourfolds containing a plane those that are known to be rational form a dense subset.
\smallskip

One day the rationality question for cubic fourfolds will be decided, in the sense of any of the above conjectures or otherwise. However, I believe that the interest in the Hodge structure $H^4(X,\ZZ)$ and in the K3 category $\ka_X$ of cubic fourfolds will persist.
For example, one might think of controlling birationality of cubic fourfolds by the
following related guess.

\begin{conj}
Assume that for two smooth cubic fourfolds $X,X'\subset\PP^5$ 
there exists an exact linear equivalence $\ka_X\cong\ka_{X'}$. Then $X$ and
$X'$ are birational. 
\end{conj}

In the spirit of \cite{AT}, this conjecture should be equivalent to  predicting birationality of $X$ and $X'$ whenever there exists an oriented Hodge isometry
$\widetilde H(X,\ZZ)\cong\widetilde H(X',\ZZ)$, see the next section for further explanations. Some evidence for the conjecture is provided by Fan and Lai \cite{FL}. The converse of the conjecture does certainly not hold.

\section{Two famous examples}

Two types of smooth cubic fourfolds have always been served as a testing ground
for various questions and conjectures: 
Pfaffian cubic fourfolds and cubic fourfolds containing a plane. The latter play a central role in \cite{AT}.

\subsection{Pfaffian cubic fourfolds} A cubic fourfold $X$ is called Pfaffian 
if it is obtained as the intersection of a linear $\PP^5\subset\PP(\bigwedge^2W^\ast)$
with the universal Pfaffian in $\PP(\bigwedge^2W^\ast)$ of all degenerate two-forms on a six-dimensional vector space $W$. It turns out that a (generic) Pfaffian cubic
contains certain quartic normal scrolls $\Sigma_P$ parametrised by points 
$P\in S\subset {\text G}(2,W)$ in a K3 surface $S$ of degree $14$. The presence
of these surfaces implies the existence of a non-trivial Hodge
class $0\ne v\in H^{2,2}(X,\ZZ)_{\rm pr}$, namely $v=4h^2-3[\Sigma_P]$. Its orthogonal complement
in $H^4(X,\ZZ)_{\rm pr}$ is Hodge isometric to the primitive cohomology of $S$:
\begin{equation}\label{eqn:PfaffianHodge}
v^\perp\cong H^2(S,\ZZ)_{\rm pr}.
\end{equation}
These Hodge theoretic computations go back to Beauville and Donagi \cite{BD}.
They also showed birationality of the Fano variety $F(X)$ of lines  and
of the Hilbert scheme $S^{[2]}$   and proved rationality of $X$,
a classical result known to Fano already.

As alluded to already in Section \ref{sec:Hodge}, from the categorical point of view,
the extension of (\ref{eqn:PfaffianHodge}) to a Hodge isometry
$$\widetilde H(X,\ZZ)\cong \widetilde H(S,\ZZ)$$
is more conceptual, but it does not fully reflect the richness of the geometric situation.\smallskip

On the categorical side, Kuznetsov \cite{Kuz5} used his homological projective duality
to construct an exact linear equivalence
$$\ka_X\cong \Db(S).$$
An alternative and geometrically more direct argument using the family of quartic scrolls $\Sigma_P$  was worked out by Addington and Lehn \cite{AL}.

\subsection{Cubics with a plane}
A smooth cubic fourfold containing a plane $\PP^2\cong P\subset X\subset\PP^5$ also admits an additional Hodge class $0\ne v\in H^{2,2}(X,\ZZ)_{\rm pr}$, namely $v=h^2-3[P]$. Its orthogonal complement
$v^\perp\subset H^4(X,\ZZ)_{\rm pr}$ is linked to the K3 surface $S$
parametrising all residual quadrics $Q_x$, $x\in S$, of linear intersections
containing the plane  $P\subset X\cap \PP^3$ together with a ruling. The rulings
glue to a Brauer--Severi variety $F_P\to S$ to which one naturally associates
a Brauer class $\alpha\in \Br(S)$.

It turns out that the Hodge structure $v^\perp$ isometrically embeds into
$H^2(S,\ZZ)_{\rm pr}$:
\begin{equation}\label{eqn:planeHodge}
v^\perp\,\,\hookrightarrow H^2(S,\ZZ)_{\rm pr}.
\end{equation}
Unlike (\ref{eqn:PfaffianHodge}), this embedding is of index two,
which was already observed by  Voisin \cite{Voisin}. However, a clear understanding in terms of twisted K3 surfaces
became only possible after \cite{HS}. Maybe the best way of viewing (\ref{eqn:planeHodge}) is either as a Hodge isometry with the twisted transcendental lattice $$v^\perp\cong T(S,\alpha)\coloneqq \ker(\alpha\colon T(S)\to \QQ/\ZZ),$$ 
which works for a very general cubic containing a plane $P\subset X$, or as a Hodge isometry
$$\widetilde H(X,\ZZ)\cong \widetilde H(S,\alpha,\ZZ).$$
We refer to \cite{HS} for the definition of the Hodge structure  of a twisted
K3 surface $(S,\alpha)$ which requires the choice of an additional $B$-field lift of $\alpha$. \smallskip

The Hodge theoretic picture is complemented by Kuznetsov \cite{Kuz1} proving the existence of an exact linear equivalence $$\ka_X\cong \Db(S,\alpha).$$
Here, the right hand side denotes the bounded derived category of $\alpha$-twisted
coherent sheaves on $S$, a notion that has its precursor in \cite{MukaiVB} and was further studied in \cite{HS}. A more geometric proof, following the idea
of Addington and Lehn \cite{AT} for Pfaffian cubics, can be found in \cite[Ch.\ 7.3]{HuyCubic}. Using the Hodge structure of the twisted K3 surface $(S,\alpha)$,
Kuznetsov \cite{Kuz1} showed that for the very general cubic containing a plane the category $\ka_X$ is not equivalent to $\Db(S')$ of any K3 surface $S'$.

\section{Hassett divisors}
Consider the moduli space $\kc$ of all smooth cubic fourfolds. It is a Deligne--Mumford stack with a 20-dimensional quasi-projective coarse moduli space.
A standard Hodge theory argument reveals that 
the locus
$$\{ X\mid H^{2,2}(X,\ZZ)_{\rm pr}\ne0\} \subset \kc$$
of all smooth cubic fourfolds admitting a non-trivial primitive Hodge class is a countable union $\bigcup \kc_d$ of irreducible divisors $\kc_d\subset\kc$. For
example, with the appropriate convention, $\kc_8$ is the set of all cubics containing a plane and $\kc_{14}$ is the set of all Pfaffian cubics.\smallskip

Hassett \cite{HassettComp} provides a detailed numerical analysis of this countable union. Firstly, the index $d$ is the discriminant of $H^{2,2}(X,\ZZ)$ of the very general $X\in \kc_d$.\footnote{If this is taken as the definition of the divisor $\kc_d$, its irreducibility is highly non-trivial.}
Then he introduces two numerical conditions ($\ast$) and ($\ast\ast$), which we will
not spell out here, such that\smallskip

$\bullet$ The divisor $\kc_d$ is non-empty
if and only if $d$ satisfies ($\ast$).\smallskip

$\bullet$ The divisor $\kc_d$ with $d$ satisfying ($\ast\ast$) 
describes the set of all cubic fourfolds $X$ for which there exists a primitive
isometric embedding of Hodge structures $H^2(S,\ZZ)_{\rm pr}\,\hookrightarrow H^4(X,\ZZ)_{\rm pr}$ (automatically of corank one), where $(S,H)$ is a polarised
K3 surface of degree $d$. From a more categorical perspective, suppressing the
polarisation, one would rephrase this as
\begin{equation}\label{eqnstarstar}
\bigcup_{{\rm (\ast\ast)}}\kc_d=\{X\mid \widetilde H(X,\ZZ)\cong \widetilde H(S,\ZZ)\}.
\end{equation}
The divisor $\kc_{14}$ of Pfaffian cubic fourfolds is one of these divisors.
Moreover, assuming Conjecture \ref{conj:Hassett}, the union (\ref{eqnstarstar})
is the set of all rational smooth cubic  fourfolds.

\begin{remark}\label{rem:Ratlatest}
Spelling out the numerical condition ($\ast\ast$), one finds that
Conjecture \ref{conj:Hassett} predicts rationality of all cubics contained in the Hassett divisors $\kc_{14}$, $\kc_{26}$, $\kc_{38}$, $\kc_{42}$, $\kc_{62}$, $\kc_{74}$,
\ldots. For $d=14$ this is classical and for the next two values $d=26$ and $d=38$ it was confirmed recently by Russo and Staglian\`o \cite{SR}. The
case $d=42$ was settled again by Russo and Staglian\`o \cite{SR2}, so that the  first Hassett divisor in this series for which rationality is not yet known is $\kc_{62}$.
\end{remark}

$\bullet$
Building upon Hassett's work, \cite{HuyComp} gives a numerical description of
all cubic fourfolds associated with twisted K3 surfaces. The numerical condition
was denoted ($\ast\ast'$). So
$$\bigcup_{{\rm (\ast\ast')}}\kc_d=\{X\mid \widetilde H(X,\ZZ)\cong \widetilde H(S,\alpha,\ZZ)\}.$$\smallskip

$\bullet$ Another numerical condition ($\ast\!\ast\!\ast$), describing the set of cubic fourfolds in Conjecture \ref{conj:GS}, was described by Addington \cite{Add}, improving earlier work of Hassett \cite{HassettComp}. \smallskip

These four numerical conditions are linked by a string of  implications
$${\rm (\ast\!\ast\!\ast) \Rightarrow (\ast\ast) \Rightarrow (\ast\ast') \Rightarrow (\ast),}$$ which are  all strict. For example, $d=74$ satisfies ($\ast\ast$) but not ($\ast\!\ast\!\ast$) and $d=8,32$ satisfy ($\ast\ast'$) but not ($\ast\ast$).\footnote{Results of Ouchi \cite{Ouchi} link ${\rm (\ast\!\ast\!\ast)}$
and $(\ast\ast')$ for certain values of $d$. For example, for $d=8$ he shows
that whenever $\alpha=1$ and hence $X$ is rational, then $F(X)$ is birational to $S^{[2]}$.}

Note that if Conjecture \ref{conj:Hassett} or, equivalently, Conjecture \ref{conj:Kuznetsov} turns out to be true, then we still have a geometrical interpretation
of the numerical condition ($\ast\!\ast\!\ast$). But which distinguished
geometric property of a cubic fourfold would then correspond to ($\ast\ast'$)?

For the reader's convenience, here are the first values of $d$ satisfying the respective four numerical conditions:\\
$$\begin{array}{ll}
{\rm(\ast)} &d=8,12,14,18,20,24,26,30,32,36,38,42,44,48,50,54,56,60,62,66,68,72,74,78, \ldots.\\[4pt]
{\rm(\ast\ast')} &d=8,14,18,24,26,32,38,42,50,54,62,72,74,78,86,96,98,104\ldots.\\[4pt]
{\rm (\ast\ast)}&d=14,26,38,42,62,74,78,86,98,114,122,134, 146,\ldots.\\[4pt]
{\rm (\ast\!\ast\!\ast)}&d=14,26,38,42,62,86,114,122,134, 146,\ldots.
\end{array}$$

\section{Addington \& Thomas and twists}

The approach of Addington and Thomas \cite{AT} involves two main ingredients.\smallskip

(i) On the one hand, surprisingly subtle arguments in lattice theory are developed to pass
from Hassett's point of departure of primitive  isometric embeddings $H^2(S,\ZZ)_{\rm pr}\hookrightarrow H^4(X,\ZZ)_{\rm pr}$ of Hodge structures
 to the more category adapted Hodge isometries
$\widetilde H(X,\ZZ)\cong\widetilde H(S,\ZZ)$. Also, it is proved that for every $d$ satisfying ($\ast$), i.e.\ such that $\kc_d$ is not empty, there exists a cubic fourfold
in the intersection $X\in \kc_8\cap \kc_d$, i.e.\ a cubic with a plane, for which moreover the associated Brauer class $\alpha$ on
the associated K3 surface is trivial.
Hence, by results of Kuznetsov \cite{Kuz1},  every non-empty Hassett divisor $\kc_d$ contains a cubic $X\in\kc_d$ 
with $\ka_X\cong\Db(S)$ for some K3 surface $S$.\smallskip

(ii) The deformation theory of the Fourier--Mukai kernel $\ke$ of a given equivalence
$\ka_X\cong \Db(S)$ is controlled by its action on cohomology. As long
as no obstruction is cohomologically detected, the FM kernel deforms. This allows
one to deform $\ke$ over an open subset $\emptyset\ne U\subset \kc_d$ for
every $d$ satisfying ($\ast\ast$). This kind of deformation argument was
first exploited by Toda \cite{Toda} and subsequently by Huybrechts--Macr\`i--Stellari \cite{HMS} and Huybrechts--Thomas \cite{HT}. The passage from formal deformations
to deformations over Zariski open sets relies on results of Lieblich \cite{Lieblich}.

As the property for a FM kernel to define an equivalence is an open condition, this leads to the main
result of \cite{AT} which was in \cite{HuyComp} extended to the twisted case following the same ideas.

\begin{thm}[Addington--Thomas, Huybrechts]
Let $X\subset \PP^5$ be a smooth cubic fourfold and $(S,\alpha)$
a twisted K3 surface. Then there exists a Hodge isometry
$\widetilde H(X,\ZZ)\cong \widetilde H(S,\alpha,\ZZ)$ if and only
if there exists an exact equivalence $\ka_X\cong \Db(S,\alpha)$.\TBC{Orientation business?}
\end{thm}

\begin{remark}
The result was not quite proved in \cite{AT,HuyComp} as stated above. Due to the limitation of the techniques, it was only proved for Zariski dense open subsets of each Hassett divisor $\kc_d$ and assuming that all equivalences are of Fourier--Mukai type.
The recent results by Bayer et al \cite{Bayer} and Li, Pertusi, and Zhao \cite{LPZ}
allow us to state the result in this more complete form.
\end{remark}

Extending results about K3 surfaces to the twisted setting might not seem very exciting or particularly difficult. There is however one main advantage of working
with general twisted K3 surfaces $(S,\alpha)$: According to \cite{HMS2}, the group $\Aut(\Db(S,\alpha))$
of exact linear auto-equivalences is often (enough) computable. Before giving a precise statement, recall that for a projective K3 surface $S$ the group $\Aut(\Db(S))$ of exact linear auto-equivalences is conjecturally described by Bridgeland's conjecture
\cite{BrK3}, which is open beyond the Picard rank one case, and expected to be always large and complicated.

\begin{thm}[Huybrechts--Macr\`i--Stellari]
Assume $(S,\alpha)$ is a twisted K3 surface such that $\Db(S,\alpha)$ does not contain any spherical objects. Then  the kernel of the natural representation
$$\Aut(\Db(S,\alpha))\to \Aut(\widetilde H(S,\alpha,\ZZ))$$ is spanned by the double shift $[2]$.
\end{thm}

It turns out that the situation occurs frequently for twisted K3 surfaces
associated with cubic fourfolds.

\begin{thm}[Huybrechts]
For an infinite number of $d$ satisfying {\rm ($\ast\ast'$)} the twisted K3 surface
$(S,\alpha)$ associated with the very general cubic $X\in \kc_d$ does not have any spherical objects.
\end{thm}

The deformation theory developed in \cite{AT,HuyComp} then allows one
to use this to describe $\Aut(\ka_X)$ for very general cubic fourfolds. Ideally,
the next result should cover all $X$ not contained in any Hassett divisor $\kc_d$,
but this is not known presently.

\begin{thm}[Huybrechts]\label{thmHuyAut}
The group of symplectic auto-equivalences $\Aut_s(\ka_X)$ of the very general
cubic $X\subset\PP^5$ contains the group of even shifts $\ZZ\cdot[2]$ as a subgroup of index three
$$\Aut_s(\ka_X)/\ZZ\cdot[2]\cong\ZZ/3\ZZ.$$
\end{thm}

\begin{remark}
There are some similarities but also differences between the K3 category
$\ka_X$ of a very general cubic $X$ and the bounded derived category $\Db(S)$
of the very general non-projective K3 surface $S$. 
For the former the Hodge classes form 
the positive-definite rank two lattice $\widetilde H^{2,2}(X,\ZZ)\cong A_2$
and for the latter  $\widetilde H^{1,1}(S,\ZZ)$ is the indefinite hyperbolic plane.

In both cases, 
the group of auto-equivalences $\Aut(\ka_X)$ and $\Aut(\Db(S))$ can be computed and
Bridgeland's conjecture is known. However, we do not have this result for all
cubics $X$ with $\widetilde H^{2,2}(X,\ZZ)\cong A_2$ but only for an unspecified very general subset of them.
\end{remark}

\section{What is left open?}

To establish in its full glory the analogue of the picture for K3 surfaces developed by Mukai \cite{MukaiVB} and Orlov \cite{OrlovK3} (and Huybrechts--Stellari \cite{HS}
in the twisted setting),
one would need to prove the following.

\begin{conj}\label{conj:main}
There exists an exact linear equivalence $\ka_X\cong\ka_{X'}$ between the K3 
categories $\ka_X$ and $\ka_{X'}$ of two smooth cubic fourfolds $X,X'\subset \PP^5$ 
if and only if there exists an orientation preserving Hodge isometry $\widetilde H(X,\ZZ)\cong \widetilde H(X',\ZZ)$.
\end{conj}

The `only if' direction, at least for equivalences of FM type, is not difficult and was proved in \cite{HuyComp}. For the `if' part, which amounts to producing an equivalence from a given
Hodge isometry, we have two types of results:\smallskip

(i) The conjecture holds for $X\in \kc_d$ with $d$ satisfying ($\ast\ast'$).
This is proved generically in \cite{AT} for $d$ satisfying the stronger
condition ($\ast\ast$) and then twisted in \cite{HuyComp} to also cover
the weaker condition ($\ast\ast'$). Again, to cover all of $\kc_d$ and general
equivalences one needs the more recent work \cite{Bayer,LPZ}.
The conjecture also holds for the very general (but unspecified) cubic in every $\kc_d$ (without any further condition on $d$).\smallskip

(ii) It was observed in \cite{HuyComp} that if $X$ is not contained in any Hassett divisor, so $X\not\in \bigcup_{\rm(\ast)}\kc_d$, any other cubic $X'$ with
$\ka_X\cong\ka_{X'}$ is actually isomorphic to $X$:
$$\ka_X\cong\ka_{X'}~\Leftrightarrow~X\cong X'.$$
Indeed, by the easy `only if' direction in Conjecture \ref{conj:main}
any such (FM) equivalence induces a Hodge isometry
between the transcendental lattices $H^4(X,\ZZ)_{\rm pr}\cong H^4(X',\ZZ)_{\rm pr}$.
By the global Torelli theorem, this implies $X\cong X'$.

Alternatively, any equivalence
$\ka_X\cong \ka_{X'}$ for $X\not\in \bigcup_{\rm(\ast)}\kc_d$ leads to an isomorphism between the  Fano varieties $F(X)\cong F(X')$ and then a geometric version of the global Torelli theorem due to Charles, cf.\ \cite[Ch.\ 2.3]{HuyCubic}, is enough to conclude.\smallskip

There are further results confirming the analogy between the Mukai--Orlov theory
for $\Db(S)$ or more generally for $\Db(S,\alpha)$ of (twisted) K3 surfaces and the theory of K3
categories $\ka_X$ of cubic fourfolds. For example, the number of cubic fourfolds  with equivalent K3 categories $\ka_X$ is always finite. This was proved in \cite{HuyComp} but also independently observed by Perry. A more precise count for very general cubic fourfolds in Hassett divisors was discussed by Pertusi \cite{Pertusi} and for certain cubics not linked to K3 surfaces by  Fan and Lai \cite{FL}.

\section{Degree shift}

We have seen in Theorem \ref{thmHuyAut} that $\Aut_s(\ka_X)/\ZZ\cdot[2]\cong\ZZ/3\ZZ$ for very general cubic fourfolds. The reason for the
appearance of the order three group is the degree shift functor. This is a certain
auto-equivalence $T_X\in\Aut(\ka_X)$ satisfying
$$T_X^3\cong [2].$$ There are two interpretations of $T_X$.
Following Kuznetsov \cite{Kuz4}, $T_X$ is the functor $$T_X\colon E\mapsto j^\ast(E\otimes\ko(1)),$$
where $j^\ast\colon\Db(X)\to\ka_X$ is the projection, i.e.\ the left adjoint of the inclusion $\ka_X\,\hookrightarrow \Db(X)$. Alternatively,
using Orlov's interpretation of $\ka_X$ as the category of graded matrix factorisations
$(K\stackrel{\alpha}{\to}L\stackrel{\beta}{\to}K(3))$, $T_X$ can be thought of
as the degree shift\footnote{One would expect that the interpretation of $\ka_X$ as a category of matrix factorisations would lead to concrete insights but so far a really convincing application to the circle of ideas discussed here is still missing.}
 $$T_X\colon (K\stackrel{\alpha}{\to}L\stackrel{\beta}{\to}K(3))\mapsto
(L\stackrel{-\beta}{\to}K(3)\stackrel{-\alpha}{\to}L(3)).$$
\begin{remark}
Note that an arbitrary K3 surface $S$ does not admit an auto-equivalence $T_S\in
\Aut(\Db(S))$ with the above property $T^3\cong [2]$. Is this property maybe characterising
those K3 surfaces $S$ for which there exists a cubic fourfold $X$ with
$\Db(S)\cong\ka_X$? Is maybe the cohomological existence of such an equivalence $T$ enough?
\end{remark}

The degree shift functor $T_X$ does not only play a special role for very general
cubic fourfolds $X$, for which $T_X$ and $[2]$ essentially generate the group of
auto-equivalences of $\ka_X$, but in fact for all. It allows one to reconstruct the cubic from $\ka_X$ via the Jacobi ring.

\begin{thm}[Huybrechts--Rennemo] Any exact linear equivalence $\ka_X\cong\ka_{X'}$ between the
K3 categories of two smooth cubic fourfolds that commutes with the degree shift functors $T_X$ and $T_{X'}$ induces an isomorphism of the graded  Jacobian rings
$J(X)\cong J(X')$ and, therefore, an isomorphism $X\cong X'$.
\end{thm}

In \cite{HR} the result was then combined to give a new proof of the global Torelli theorem for cubic fourfolds. One first uses Conjecture \ref{conj:main} for the very general cubic $X$ to lift any Hodge isometry $H^4(X,\ZZ)_{\rm pr}\cong H^4(X',\ZZ)_{\rm pr}$ to an equivalence $\ka_X\cong\ka_{X'}$ and then Theorem \ref{thmHuyAut} to show that such an equivalence automatically commutes
with the degree shift functors $T_X$ and $T_{X'}$. This proves a generic global Torelli
theorem which for cubic fourfolds implies the global Torelli theorem for all
smooth cubic fourfolds. Another proof of the theorem for very general cubics but 
assuming only a cohomological compatibility with the degree shift functor
was given in \cite{BLMS}.

\TBC{Brakkee,.}

\section{The structure of $\ka_X$}

The K3 category $\ka_X$ gives rise to interesting higher-dimensional hyperk\"ahler categories $\ka_X^{[n]}$. They are constructed as ${\mathfrak S}_n$-equivariant
versions of the exterior product $$\ka_X^n=\ka_X\boxtimes\cdots\boxtimes\ka_X\subset
\Db(X)\boxtimes\cdots\boxtimes\Db(X)\cong\Db(X^n)$$ or, more directly, as a certain full subcategory of $\Db([X/{\mathfrak S_n}])$.

Consider  a projective moduli space $M$ of objects in $\ka_X$, e.g.\ say stable with respect to some stability condition, and assume $M$ is of dimension $2n$.

\begin{question} Is $\ka_X^{[n]}$ equivalent to $\Db(M)$? 
\end{question}

By virtue of results of Li, Pertusi, and Zhao  \cite{LPZ1}, the Fano variety of lines $F(X)$ can be viewed as such a moduli space. Thus, as a special case of this general question one
recovers a conjecture attributed to Galkin who asked: Is there an
exact linear equivalence $$\ka_X^{[2]}\cong\Db(F(X))?$$

Note that for the general cubic $X\in\kc_d$ with $d$ satisfying ($\ast\!\ast\!\ast$) one has
$F(X)\cong S^{[2]}$ and $\ka_X\cong\Db(S)$ which together indeed give
$$\ka_X^{[2]}\cong\Db(S)^{[2]}\cong\Db(S^{[2]})\cong \Db(F(X)).$$

Certain evidence for an affirmative answer for $n=2$ is 
the equality in the Grothendieck group of triangulated categories
$$[\ka_X^{[2]}]=[\Db(F(X))]\in K_0(\text{dg-cat}),$$
which is rather easy to prove. A stronger result  is due to
 Belmans, Fu, and Raedschelders \cite{BFR} showing that 
 $\ka_X^{[2]}$ and $\Db(F(X))$ can both be realised as semi-orthogonal factors
 of two semi-orthogonal decompositions of $\Db(X^{[2]})$ with all other factors individually equivalent to $\Db(X)$.
 
 \begin{remark}
 Besides being very suggestive, the above question would also shed light on 
 the following open question for K3 surfaces. Assume $S$ is a K3 surface and $M(v)$
 is a projective moduli space of stable sheaves. Is its derived category
 $\Db(M(v))$ equivalent to the derived category $\Db(S^{[n]})$ of the Hilbert scheme
 of the same dimension? In other words, are all projective moduli spaces of stable
 sheaves of the same dimension on a fixed K3 surface $S$ derived equivalent to each other? Naive attempts to settle this question either way have all failed. However, it has been observed in concrete examples
 \cite{Add1,Add2} and  Beckmann and Bottini could show that Taelman's Mukai lattices
 of $M(v)$ and $S^{[n]}$ are Hodge isometric.
 \end{remark}
 
\section{Unknown deformations}

It seems that at least conjecturally the picture is quite clear. However, there is
a big part of the deformation theory of the cubic that seems unaccounted for.

Let us begin with the case of K3 surfaces. The classical first-order deformations of $S$ as a complex K3 surface are parametrised by $H^1(S,\kt_S)\cong\CC^{20}$.
If a polarisation of $S$ is fixed, those deformations that preserve the polarisation
form a hyperplane $$\CC^{19}\cong H^1(S,\kt_S)_{\rm pol}\subset H^1(S,\kt_S)\cong\CC^{20}.$$ The non-commutative deformations
of $S$, of which we think as deformation of its derived category $\Db(S)$, are
parametrised by the Hochschild cohomology $H\!H^2(\Db(S))$. The identification
\begin{equation}\label{eqn:HH}
H\!H^2(\Db(S))\cong H^2(S,\ko)\oplus
H^1(S,\kt_S)\oplus H^0(S,{\bigwedge}^2\kt_S)
\end{equation}
reveals that this is a space of dimension $22$. The one-dimensional subspaces
$H^2(S,\ko)$ and $H^0(S,\bigwedge^2\kt_S)$ correspond to 
twisted and Poisson deformations of $S$, respectively. Note, however, that the direct
sum decomposition is typically not preserved by $\Aut(\Db(S))$.

Assume now that the K3 surface $S$ is associated with a cubic fourfold $X$, i.e.\ $\Db(S)\cong \ka_X$. Then the induced map $H^1(X,\kt_X)\to H\!H^2(\Db(S))$
is injective and its image contains the hyperplane of polarised deformations
$$H^1(S,\kt_S)_{\rm pol}\subset H^1(X,\kt_X)\,\,\hookrightarrow H\!H^2(\Db(S)).$$ The one remaining deformation direction of $X$, i.e.\ the one not contained in $H^1(S,\kt_S)_{\rm pol}$, will usually involve all three summands in (\ref{eqn:HH}).

There are however non-commutative dimensions of $X$ as well. It turns out
that\footnote{I wish to thank Pieter Belmans for this information and an intersting email conversation.}
$$\begin{array}{ccccccc}H\!H^2(X)&\cong& H^2(X,\ko)&\oplus& H^1(X,\kt_X)&\oplus &H^0(X,{\bigwedge}^2\kt_X)\\[5pt]
&\cong &0&\oplus& \CC^{20}&\oplus& \CC^{20}.
\end{array}$$
The curious numerical coincidence of both non-trivial subspaces being of dimension $20$ is unexplained.

Hodge theory suggests that one more non-commutative deformation of $\Db(S)$
is accounted for by non-commutative deformations of $X$. This leaves us with
one non-commutative deformation of the K3 surface $S$ and $19$ (infinitesimal) non-commutative
deformations of the cubic fourfold $X$ that seem unrelated.


\end{document}